\documentclass[10pt]{amsart}
\usepackage{amsmath,amstext,amssymb,amsopn,amsthm,mathrsfs}
\usepackage{graphicx}
\usepackage[english]{babel}
\usepackage{times}
\usepackage{fancybox}
\usepackage{epsfig}
\usepackage{graphics}
\usepackage{color}
\usepackage{amsmath}
\usepackage{amsfonts}
\usepackage[ansinew]{inputenc}
\usepackage{tikz}

\numberwithin{equation}{section}

\begin{document}

\title[An non-centered asymmetric Cantor-like Set] {An non-centered asymmetric Cantor-like Set}

\author{Lauren Wszolek}
\address{Department of Mathematical and Actuarial Sciences, Roosevelt University  Chicago, Il, 60605, USA.}
\email{[Lauren Wszolek] lwszolek@mail.roosevelt.edu}
\author{Wilfredo O. Urbina}
\email{[Wilfredo Urbina]wurbinaromero@roosevelt.edu}

\begin{abstract}
The ternary Cantor set $C$, constructed by George Cantor in 1883, is probably the best known example of a perfect nowhere-dense set in the real line, but as we will see later, it is not the only one. The present article will delve into the richness and the peculiarities of $C$ through exploration of several variants and generalizations, and will provide an example of a non-centered asymmetric Cantor-like set.
\end{abstract}
\maketitle

\section{Introduction: Cantor Ternary Set.}
The (ternary) Cantor set $C$ is probably the best known example of a {\em perfect nowhere-dense} set in the real line. It was constructed by George Cantor in 1883, see \cite{cantor}.\\

$C$ is obtained from the closed interval $[0,1]$ by a sequence of deletions of open intervals known as ''middle thirds".
We begin with the interval $[0,1]$, let us call it $C_0$, and remove the middle third, leaving us with leaving us with the union of two closed intervals of length $1/3$,
$$C_1= \left[0,\frac{1}{3}\right]  \cup \left[\frac{2}{3}, 1\right] .$$ Next we remove the middle third from each of these intervals, leaving us with the union of four closed intervals of length $1/9$,
\begin{equation*}
C_2=\left[0,\frac{1}{9}\right] \cup \left[\frac{2}{9},\frac{1}{3}\right] \cup\left[\frac{2}{3},\frac{7}{9}\right] \cup\left[\frac{8}{9},1\right] .
\end{equation*}
Then we remove the middle third of each of these intervals leaving us with eight intervals of length $1/27$,
\begin{equation*}
C_{3}= [0,\frac{1}{27}]\cup[\frac{2}{27}, \frac{1}{9}]\cup[\frac{2}{9}, \frac{7}{27}]\cup [ \frac{8}{27},\frac{1}{3}] \cup [\frac{2}{3},\frac{19}{27}]\cup[\frac{20}{27},\frac{7}{9}]\cup[\frac{8}{9},\frac{25}{27}]\cup[\frac{26}{27},1].
\end{equation*}

We continue this process inductively. Then for each $n=1,2, 3\cdots $, we get a set $C_n$  which is the union of $2^n$ closed intervals of length $1/3^n$. Finally,  the {\em Cantor ternary set} $C$ is defined as the intersection
\begin{equation}
C= \bigcap_{n=0}^\infty C_n.
\end{equation}

Clearly $C \neq \emptyset $, since  trivially $0,1 \in C$. Moreover, $C$ is a closed set, being the countable intersection of closed sets, and trivially bounded, since it is a subset of $[0,1]$. Therefore, by the Heine-Borel theorem, $C$ is a {\em compact set}. Moreover, observe by the construction that if $y$ is the endpoint of some closed subinterval of a given  $C_n$ then it is also the endpoint of some of the subintervals of $C_{n+1}$. At each stage, since endpoints are never removed, it follows that $y \in C_n$ for all $n$. Thus $C$ contains all the endpoints of all the intervals that make up each of the sets $C_n$ (or alternatively,  the endpoints to the intervals removed),  all of which are rational ternary numbers in $[0,1]$, i.e. numbers of the form $k/3^n$. But $C$ contains much more than that: it is in fact an uncountable set since it is a {\em perfect set}. \footnote{A perfect set $P$ is a set that is closed and every point  $x \in P$ is a limit point, i.e there is a sequence $\{x_n\} \subset P$, $x_n \neq x$ and $x_n \rightarrow x$.} To prove it is perfect, simply observe that every point of $C$ is approachable arbitrarily closely by the endpoints of the intervals removed (thus for any $x\in C$ and for each $n \in \mathbb{N}$ there is an endpoint, let us call it $y_n \in C_n$, such that $|x-y_n| < 1/3^n$).

There is an alternative characterization of $C$, the {\em ternary expansion characterization}.
Consider the ternary representation for  $x \in [0,1]:$ \footnote{Observe, for the ternary rational  numbers $k/3^n$ there are two possible ternary expansions, since
$$ \frac{k}{3^n} = \frac{k-1}{3^n} + \frac{1}{3^n}  =  \frac{k-1}{3^n} + \sum_{k=n+1}^\infty \frac{2}{3^k}.$$ 
Similarly, for the dyadic rational numbers $k/2^n$ there are two possible dyadic expansions as
$$ \frac{k}{2^n} = \frac{k-1}{2^n} + \frac{1}{2^n}  =  \frac{k-1}{2^n} + \sum_{k=n+1}^\infty \frac{1}{2^k}.$$ 
Thus for the uniqueness of the dyadic and the ternary representations we will take  the infinite expansions representations for the dyadic and ternary rational numbers.}
 
\begin{equation}\label{ternaryexp}
 x = \sum_{k=1}^\infty \frac{\varepsilon_k(x)}{3^k}, \quad \varepsilon_k(x)=0, 1, 2 \quad \mbox{for all} \,k = 1, 2, \cdots.
\end{equation}

Observe that removing the elements where at least one of the \(\varepsilon_{k} =1\)  is the same as removing the middle third in the iterative construction, and thus 
the Cantor ternary set is the set of numbers in $[0,1]$ that can be written in base 3 without using the digit 1:
\begin{equation}\label{ternaryChar}
 C=\left\{x\in [0,1] : x = \sum_{k=1}^\infty \frac{\varepsilon_k(x)}{3^k}, \quad \varepsilon_k(x)=0, 2 \quad \mbox{for all} \; k = 1, 2, \cdots
 \right\}.
\end{equation}

Using this characterization  of  $C$ we can get a direct proof that the set is uncountable. Define 
the mapping $f: C \rightarrow [0,1]$ for $  x=\sum_{k=1}^\infty \frac{\varepsilon_k(x)}{3^k} \in C$ as 
\begin{equation}\label{cantorfun}
f(x) =    \sum_{k=1}^\infty \frac{\varepsilon_k(x)/2}{2^k}= \frac{1}{2}    \sum_{k=1}^\infty \frac{\varepsilon_k(x)}{2^k}.
\end{equation}

By observing that when $\varepsilon_k =0 \text{ or } 2$, then $\varepsilon_k/2 = 0 \text{ or } 1$ respectively, it is clear that $f$ is one-to-one correspondence from $C$ to $[0,1]$. 
As we have seen before, the uncountability of $C$ can be also obtained from the fact that $C$ is perfect, see Abbott \cite{abot}, page 90. \\

$C$ is a {\em nowhere-dense} set, meaning there are no intervals included in $C$. One way to demonstrate this is when given two arbitrary points in $C$, we can always find a number between them that requires the digit 1 in its ternary representation. Therefore there are no intervals included in $C$, and thus $C$ is a nowhere-dense set. Alternatively, we can prove this simply by contradiction. Assuming that there is a interval $I=[a,b] \subset C, \; a <b$, then $I=[a,b] \subset C_n$ for all $n$, but as $|C_n| \rightarrow 0$ and $n \rightarrow \infty$, then $|I| = b-a =0.$\\

$C$ has measure zero, since its length can be obtained after subtracting from $1$ the sum of the length of all  open intervals removed in constructing it, as given by
$$ m(C) = 1 - \sum_{n=1}^\infty \frac{2^{n-1}}{3^n} = 1 -\frac{1}{3}  \sum_{n=0}^\infty (\frac{2}{3})^{n} = 1 -  \frac{1/3}{1- 2/3}=1-1=0.$$\\

Another important property of the Cantor set is its self-similarity (i.e. fractal characteristic) across scales, as illustrated in the following figure:
\begin{center}
\includegraphics[width=3in]{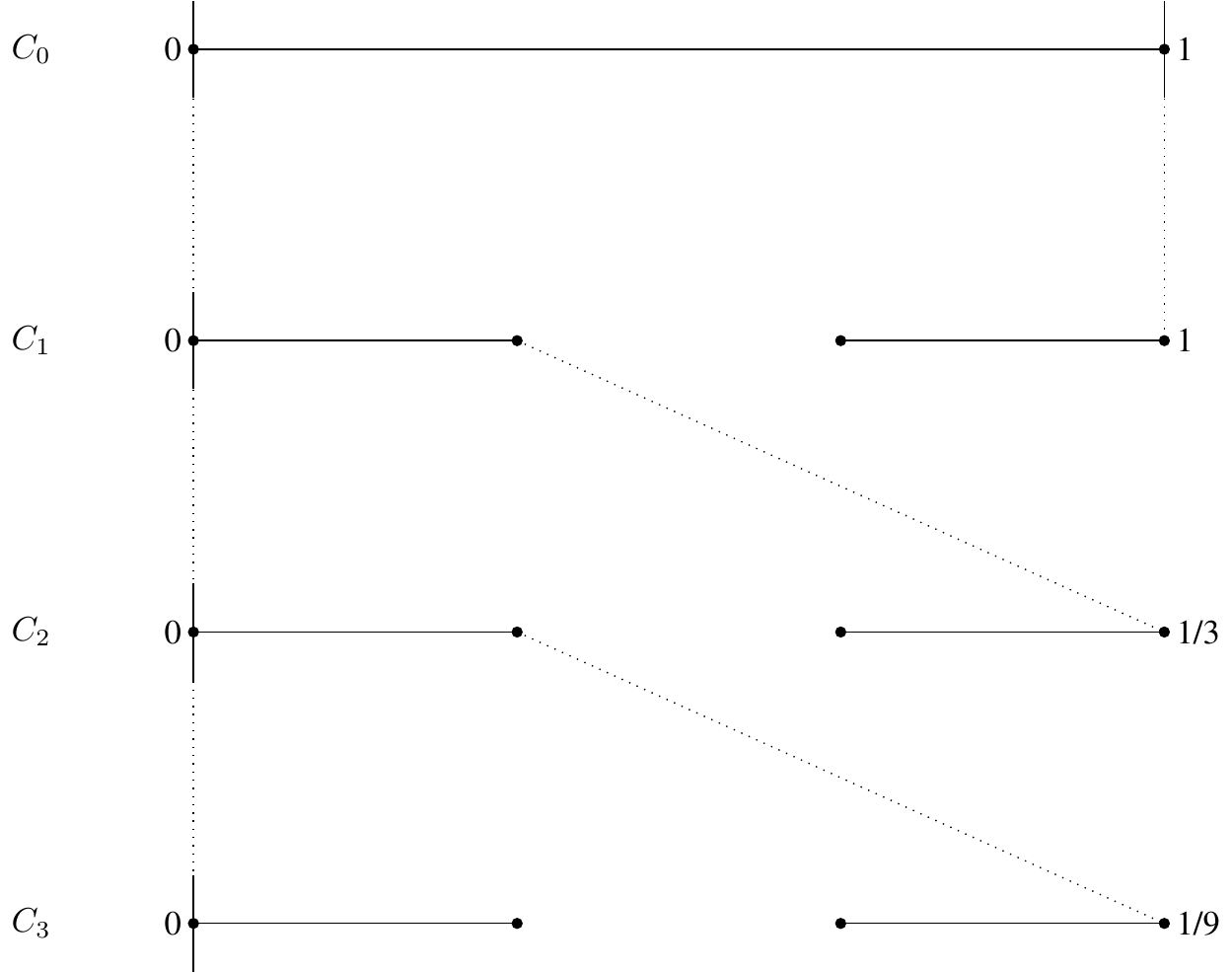}
\end{center}
 
 The Cantor ternary set $C$ was not the first perfect nowhere-dense set in the real line to be constructed. The first construction was done by the British mathematician Henry J. S. Smith in 1875, but few mathematicians were aware of Smith's construction. Vito Volterra, still a graduate student in Italy, also showed how to construct such a set  in 1881, but he published his result in an Italian journal not widely read. In 1883 Cantor rediscovered this construction himself, and due to his prestige, the Cantor ternary set became the typical example of a perfect nowhere-dense set. Following D. Bresoud \cite{bres}, we will refer to the Smith-Volterra-Cantor sets, or the $SVC(m)$ sets, as the family of examples of perfect, nowhere-dense sets, exemplified by the work of Smith, Voterra, and Cantor, where during the $n^\text{th}$-iteration, an open interval of length $1/m^n$ is removed from the center of the remaining closed intervals. Observe that $C = SVC(3)$.\\

As it was noticed in \cite{dimartinourb2}, the number $3$ is somehow special. $C$ can be obtained equivalently by removing a fixed proportion (one third) of each subinterval in each of the iterative steps, by removing the length $1/3^n$ from each subinterval in the $n^\text{th}-$step, or removing the digit ``1'' of the ternary expansion. 

For instance, if we choose $2$ instead of $3$, using the proportional construction and removing the open ``middle half" from each component, at the end of the iterative process, we get a set we will denote as $C^{1/2}$.  Below are the first three iterative steps:

\begin{eqnarray*}
C^{1/2}_{1}&=& [0,\frac{1}{4}]\cup [\frac{3}{4},1],\\
C^{1/2}_{2}&=& [0,\frac{1}{16}]\cup[\frac{3}{16}, \frac{1}{4}]\cup [\frac{3}{4},\frac{13}{16}]\cup[\frac{15}{16},1],\\
C^{1/2}_{3}&=& [0,\frac{1}{64}]\cup[\frac{3}{64}, \frac{1}{16}]\cup[\frac{3}{16}, \frac{13}{64}]\cup [ \frac{15}{64},\frac{1}{4}] \cup [\frac{3}{4},\frac{49}{64}]\cup[\frac{51}{64},\frac{13}{16}]\cup[\frac{15}{16},\frac{61}{64}]\cup[\frac{63}{64},1].
\end{eqnarray*}
Continuing this process inductively, for each $n=1,2, \cdots $, we get the set $C^{1/2}_n$ made of the union of $2^n$ closed intervals of length $1/4^n$:
$$ C^{1/2} = \bigcap_{n=1}^\infty C^{1/2}_{n}.$$
$C^{1/2}$ shares the same properties as $C$: it is a perfect, nowhere-dense set in the real line. $C^{1/2}$ also has measure zero, as its length can be obtained by subtracting from $1$ the sum of the length of all  open intervals removed in constructing it:
\begin{eqnarray*}
1- \frac{1}{2} -2 \frac{1}{8}-4\frac{1}{32} +\cdots &=& 1- \frac{1}{2}[1+ \frac{1}{2}+\frac{1}{4} +\cdots ]\\
&=& 1- \frac{1}{2}  \sum_{n=0}^\infty (\frac{1}{2})^n =  1- \frac{1}{2}  \frac{1}{1-\frac{1}{2}} =  1- \frac{1}{2} \frac{1}{\frac{1}{2}} = 1- 1=0.
\end{eqnarray*}

Also observe that we can get a {\em expansion characterization} of $C^{1/2}$, but in base $n=4$: 
\begin{equation*}
 C^{1/2}=\left\{x\in [0,1] : x = \sum_{n=1}^\infty \frac{\varepsilon_n(x)}{4^n}, \quad \varepsilon_n(x)=0, 3 \quad \mbox{for all} \; n = 1, 2, \cdots
 \right\}.
\end{equation*}

Now, what if we remove the length $1/2^n$ from the center of each subinterval in the $n^\text{th}$ step?  The first step is the same as before, obtaining the set
$$ [0,\frac{1}{4}]\cup [\frac{3}{4},1].$$
Next we need to remove two open intervals of length $1/2^2 = 1/4$, but we get only four points:
$$\{0, \frac{1}{4}, \frac{3}{4},1\}.$$
Therefore, the process stops after this step since there are no more intervals to remove from. Clearly, this is not the same as removing the open ``middle half" from each component.\\

When $n=4$, both the proportional and the power construction can be iterated an infinity number of times, but, as we are going to see, they turn out to be different sets.

\begin{enumerate}
\item [i)] {\em Proportional construction}: If we repeat the Cantor set's construction starting with the interval [0,1] and $n=4$, and removing the open ``middle fourth" from the center of the remaining closed intervals, we get at the end of the iterative process a set that we will denote as $ C^{1/4}$. Following are the first three iterations:
\begin{eqnarray*}
C^{1/4}_1 &=& [0,\frac{3}{8}] \cup [\frac{5}{8}, 1],\\
C^{1/4}_2 &=& [0,\frac{9}{64}]\cup [\frac{15}{64}, \frac{3}{8}]\cup[\frac{5}{8},\frac{49}{64}]\cup[\frac{55}{64},1],\\
C^{1/4}_3 &=& [0,\frac{27}{512}]\cup [\frac{45}{512}, \frac{9}{64}]\cup[\frac{15}{64},\frac{147}{512}]\cup[\frac{165}{512},\frac{3}{8}]\cup[\frac{5}{8},\frac{347}{512}]\cup[\frac{365}{512},\frac{49}{64}]\cup[\frac{55}{64},\frac{467}{512}]\cup[\frac{485}{512},1].
\end{eqnarray*}
Continuing this process inductively, for each $k=1,2, \cdots, $ we get a set $C^{1/4}_k$  which is the union of $2^k$ closed intervals of length $(3/8)^k$: 
$$ C^{1/4} = \bigcap_{n=1}^\infty C^{1/4}_{n}.$$
\begin{center}
\includegraphics[width=3in]{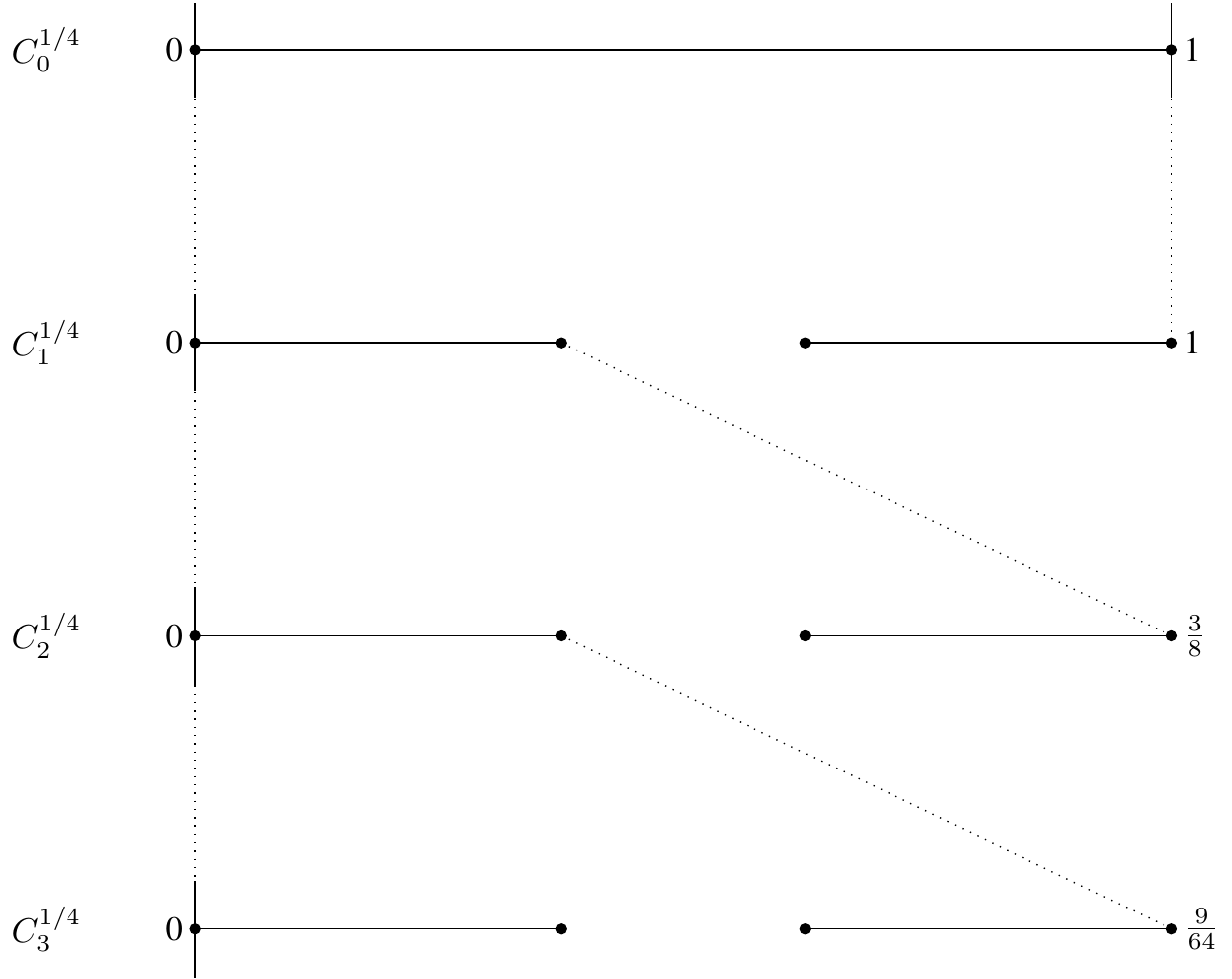}
\end{center}
$C^{1/4}$ shares the same properties as $C$ because it is a perfect, nowhere-dense set in the real line.
$C^{1/4}$ also has measure zero, as its length can be obtained after subtracting from $1$ the sum of the length of all  open intervals removed in constructing it:
\begin{eqnarray*}
1- \frac{1}{4} -2 \frac{3}{32}-4\frac{9}{256} +\cdots &=& 1- \frac{1}{4}[1+ \frac{3}{4}-\frac{9}{16} +\cdots ]\\
&=& 1- \frac{1}{4} \frac{1}{1-\frac{3}{4}} =  1- \frac{1}{4} \frac{1}{\frac{1}{4}} = 1-1 =0.
\end{eqnarray*}

Observe that $C^{1/4}$ can not be characterized using {\em expansion characterization}. It is easy to check that the ``natural base" $n=8$ does not work.\\

\item [ii)] {\em Power construction}: if we repeat the Cantor construction starting with the interval [0,1], removing in the $n$-iteration an open interval of length $1/4^n$ from the center of the remaining closed intervals, we get 
\begin{eqnarray*}
SVC(4)_1 &=& [0,\frac{3}{8}] \cup [\frac{5}{8}, 1]\\
SVC(4)_2 &=& [0,\frac{5}{32}]\cup [\frac{7}{32}, \frac{3}{8}]\cup[\frac{5}{8},\frac{25}{32}]\cup[\frac{27}{32},1]\\
SVC(4)_3 &=& [0,\frac{9}{128}]\cup [\frac{11}{128}, \frac{5}{32}]\cup[\frac{7}{32},\frac{37}{128}]\cup[\frac{39}{128},\frac{3}{8}]\cup[\frac{5}{8},\frac{89}{128}]\cup[\frac{91}{128},\frac{25}{32}]\cup[\frac{27}{32},\frac{59}{64}]\cup[\frac{119}{128},1].
\end{eqnarray*}
Then
$$ SVC(4) = \bigcap_{n=1}^\infty SVC(4) _{n}.$$
\begin{center}
\includegraphics[width=3in]{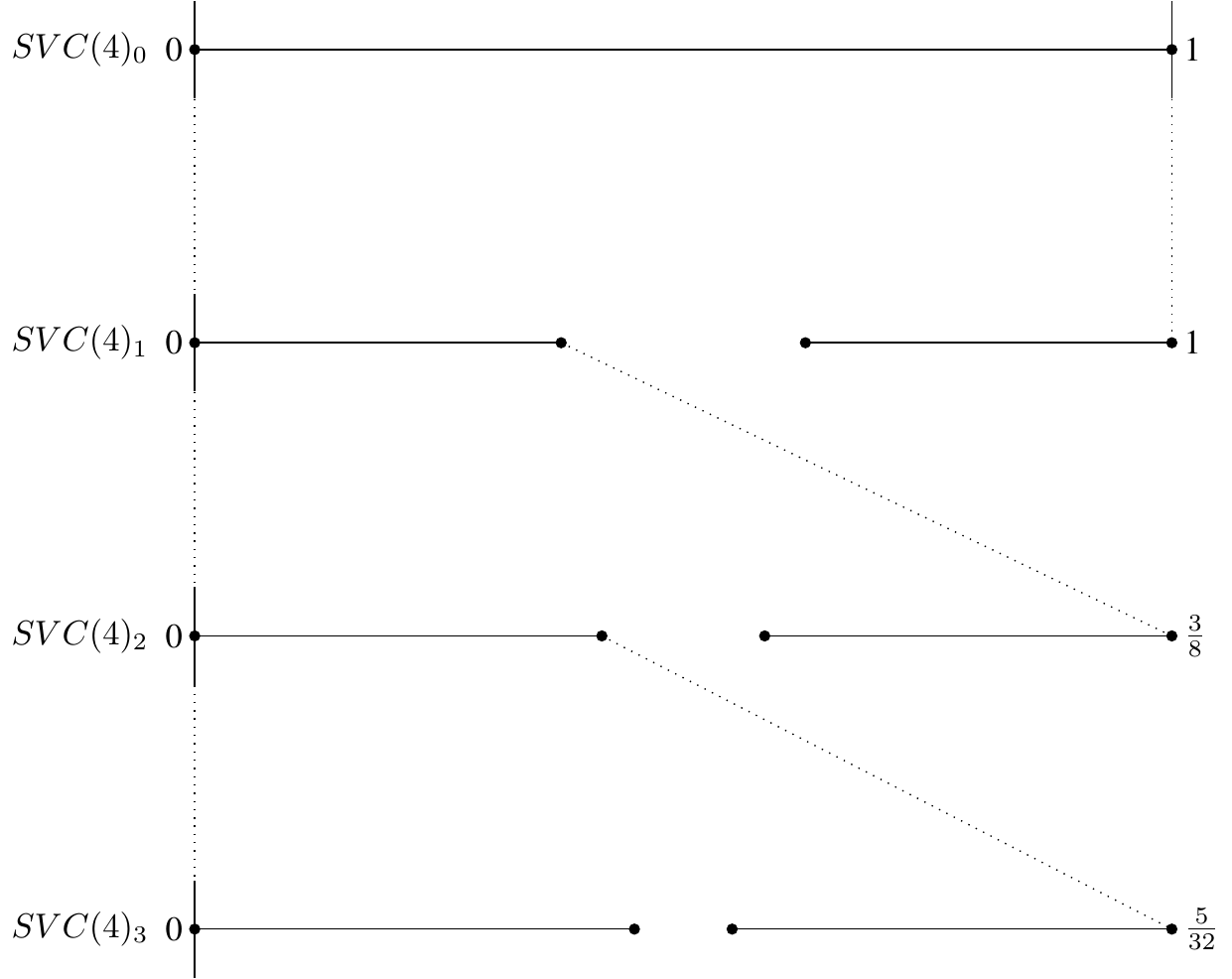}
\end{center}
 The total length of $SVC(4)$ can be obtained after subtracting from $1$ the sum of the length of all  open intervals removed in constructing it:
\begin{eqnarray*}
1- \frac{1}{4} -2 \frac{1}{16}-4\frac{1}{64} +\cdots &=& 1- \frac{1}{4}[1+ \frac{1}{2}-\frac{1}{4} +\cdots ]\\
&=& 1- \frac{1}{4} \frac{1}{1-\frac{1}{2}} =  1- \frac{1}{4} \frac{1}{\frac{1}{2}} = 1-\frac{1}{2} = \frac{1}{2}.
\end{eqnarray*}\\

Thus, the set $SVC(4)$ has positive measure equal to $1/2$. Cantor-like sets with positive measure are called {\em fat-Cantor sets}.\footnote{Observe that the only $SCV(n)$ set that has measure zero is $C$, since the total length of what is removed is
$$ 1 - \frac{1}{n} -\frac{2}{n^2} - \frac{4}{n^3}- \cdots = 1 - \frac{1}{n}[1 -\frac{2}{n} - \frac{2^2}{n^2}- \cdots ] =1 - \frac{1}{n} \frac{1}{1- 2/n} = \frac{n-3}{n-2}.$$}\\

$SVC(4)$ is called the {\em Volterra set}. This is the set that was considered in 1881 by V.  Volterra to construct his famous counter-example of a function with a bounded derivative that exists everywhere, but the derivative is not Riemann integrable in any closed bounded interval- that is, the Fundamental Theorem of Calculus fails!\\

Observe that there is not an expansion characterization of the set $SVC(4)$.\\
\end{enumerate}

If we repeat the Cantor set's {\em proportional construction} starting with the interval [0,1], removing the open ``middle three fourth" from the center of the remaining closed intervals, then we get in the first three iterations,
\begin{eqnarray*}
C^{3/4}_1 &=& [0,\frac{1}{8}] \cup [\frac{7}{8}, 1]\\
C^{3/4}_2 &=& [0,\frac{1}{64}]\cup [\frac{7}{64}, \frac{1}{8}]\cup[\frac{7}{8},\frac{57}{64}]\cup[\frac{63}{64},1]\\
C^{3/4}_3 &=& [0,\frac{1}{512}]\cup [\frac{7}{512}, \frac{1}{64}]\cup[\frac{7}{64},\frac{57}{512}]\cup[\frac{63}{512},\frac{1}{8}]\cup[\frac{7}{8},\frac{449}{512}]\cup[\frac{455}{512},\frac{57}{64}]\cup[\frac{63}{64},\frac{505}{512}]\cup[\frac{511}{512},1].
\end{eqnarray*}
Continuing this process inductively, for each $n=1,2, \cdots $ we get a set $C^{3/4}_kn$  which is the union of $2^n$ closed intervals of length $(3/8)^n$, and
$$ C^{3/4} = \bigcap_{n=1}^\infty C^{3/4}_{n}.$$
$C^{3/4}$ shares the same properties as $C$, since it is a perfect, nowhere-dense set in the real line.
Additionally, $C^{3/4}$ has measure zero, since
\begin{eqnarray*}
1- \frac{3}{4} -2 \frac{3}{32}-4\frac{3}{256} +\cdots &=& 1- \frac{3}{4}[1+ \frac{1}{4}-\frac{1}{16} +\cdots ]\\
&=& 1- \frac{3}{4} \frac{1}{1-\frac{1}{4}} =  1- \frac{3}{4} \frac{1}{\frac{3}{4}} = 1-1 =0.
\end{eqnarray*}

 Observe that in this case there is no {\em power construction.} Additionally, it is possible to create an {\em expansion characterization} of $C^{3/4}$ with base $8$:
\begin{equation}\label{Cant2}
 C^{3/4}=\left\{x\in [0,1] : x = \sum_{n=1}^\infty \frac{\varepsilon_n(x)}{8^n}, \quad \varepsilon_n(x)=0, 7 \quad \mbox{for all} \; n = 1, 2, \cdots
 \right\}.
\end{equation}
 \\

A more systematic study of the possible generalizations of the Cantor ternary set can be found in \cite{dimartinourb}.\\

\section{A non-centered asymmetric Cantor-like set}

We are going to consider a  non-centered asymmetric Cantor-like set, which we will denote $AC$, that is a hybrid of the previous examples considered. 

In the first step, we divide $[0,1]$ into four equal length intervals
$$\left[0,\frac14\right], \left[\frac14,\frac12\right],\left[\frac12,  \frac34\right],\left[\frac34, 1\right]$$ and remove the open third one, ending with
$$AC_1=\left[0,\frac12\right] \cup \left[\frac34,1\right].$$
In the second step, we divide each of the intervals in $AC_1$ into four equal length intervals, then remove the third one, and so iterate this procedure. Thus, the next two iterations:
$$\begin{aligned}  AC_{2} &=\left[0, \frac{1}{4}\right] \cup\left[\frac{3}{8}, \frac{1}{2}\right] \cup\left[\frac{3}{4}, \frac{7}{8}\right] \cup\left[\frac{15}{16}, 1\right] \\ AC_{3} &=\left[0, \frac{1}{8}\right] \cup\left[\frac{3}{16}, \frac{1}{4}\right] \cup\left[\frac{3}{8}, \frac{7}{16}\right] \cup\left[\frac{15}{32}, \frac{1}{2}\right] \cup\left[\frac{3}{4}, \frac{13}{16}\right] \cup\left[\frac{27}{32}, \frac{7}{8}\right] \cup\left[\frac{15}{16}, \frac{31}{32}\right] \cup\left[\frac{63}{64}, 1\right] . \end{aligned}$$
The graphs of them are illustrated below
\begin{center}
\includegraphics[width=6.5in]{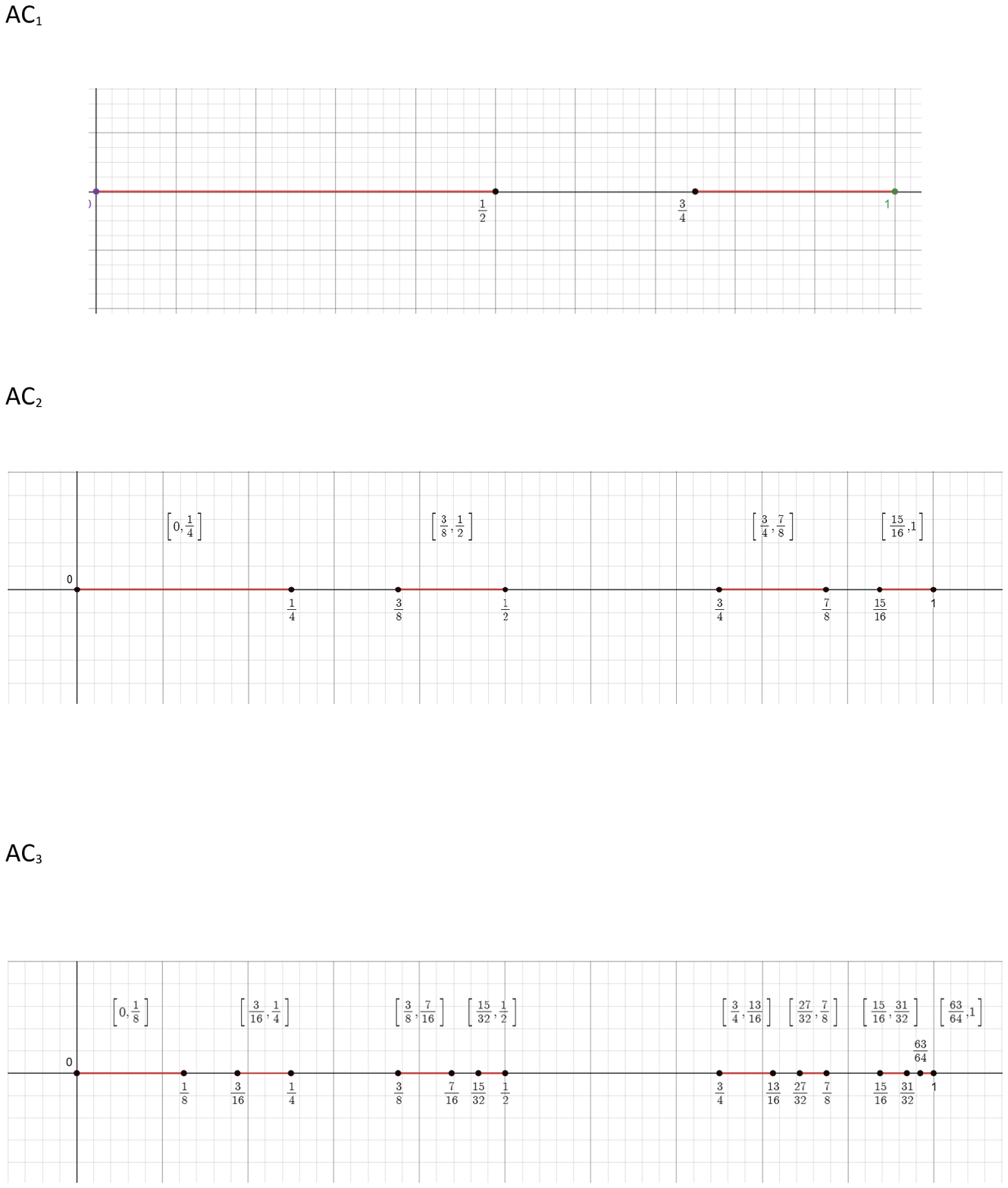}
\end{center}
Observe that $ AC_{n}$  is the union of $2^n$ closed intervals, but in contrast to the previous cases, these intervals are of differing lengths since they are of different dyadic scales, and as we are going to explain later, have total length $\frac14\cdot\left(\frac34\right)^n$.
Then we obtain the set $AC$ as the intersection of all $ AC_{n}$ 
$$AC=\bigcap_{n=0}^{\infty} AC_{n}.$$
Clearly from the construction, $AC$ is not a centered nor symmetric set, quite different from the Cantor-like sets considered above. Moreover, the set $AC$ is, in a way, a hybrid between several of the Cantor-like sets considered above.\\

The length of $AC$ can be computed by subtracting 1, the length of $[0,1]$, from the sum of the lengths of the removed intervals.
In step $1$ we remove an interval of length  $\frac{1}{4}$.  In step $2$, we remove two intervals with total length $ \frac18+\frac1{16}=\frac3{16}= \frac14 \cdot \frac3{4}$
and in step 3, we removed four intervals with total length  $\frac1{16} + \frac1{32}+\frac1{32}+\frac1{64}= \frac{9}{64}= \frac14 \cdot \frac9{16} =\frac14 \cdot \left(\frac3{4}\right)^2,$ and so on. Then, in step $n$ we remove $2^n$ intervals  with total length $\frac14\cdot\left(\frac34\right)^n$,  and therefore
$AC$ has then measure zero since
$$
1-\frac{1}{4} \sum_{n=1}^{\infty}\left(\frac{3}{4}\right)^{n}=1-\frac{\frac{1}{4}}{1-\frac{3}{4}}=1-1=0.
$$
Even though $AC$ has measure zero, again different from the previous cases, it does not have a fractal structure, i.e. it does not have a self-similar structure across the scales nor can it be characterized using any expansion characterization since the construction uses different scales in each iterative step.\\

Additionally, observe that $AC$ is nowhere-dense, since if we assume that $AC$ contains a non-trivial interval $[a, b]$ i.e. $b-a>0$, then, by the Archimedian principle, there exists a real number $n$ such that $\frac1{n}<b-a$. Then there exists $n$ such that in $AC_{n}$ there is an interval of length  $\frac1{n}$ such that $b$ and $a$ are not in the same interval, and therefore $[a,b]$ is not contained in $AC_{n}$, which is a contradiction. Thus, since $a$ and $b$ are arbitrarily chosen, $AC$ contains nontrivial intervals. Also, trivially $AC$ is a compact set and perfect since every point $ x \in AC$ is a limit point: simply take the sequence $\{x_n\}$ of extreme points of the corresponding intervals in $AC_n$ where $x$ belongs to see this.\\

A ``reflected'' set from $AC$ can be obtained if,  in the first step, after dividing the interval $[0,1]$ in four equal length intervals 
$$\left[0,\frac14\right], \left[\frac14,\frac12\right],\left[\frac12, \frac34\right],\left[\frac34, 1\right],$$  then the open second interval is removed 
$$AC_1=\left[0,\frac14\right] \cup \left[\frac12,1\right]$$  
and then the argument is iterated as previously demonstrated.\\

Finally, similar Cantor-like sets can be constructed for any $n \geq 5$. For instance, if $n=5$, in the first step we divide the interval $[0,1]$ into five equal length intervals, and then we can  remove the open fourth one or open union of the third and the fourth one, getting 
$$AC^5_1=\left[0,\frac35\right] \cup \left[\frac45,1\right],$$  
or
$$AC^5_1=\left[0,\frac25\right] \cup \left[\frac45,1\right],$$  
and then the argument can be iterated as before.

\end{document}